\documentclass[11pt]{article}
\usepackage{latexsym}
\usepackage{amssymb}
\usepackage[matrix,arrow]{xy}

\newtheorem{teo}{THEOREM}[section]
\newtheorem{prop}[teo]{PROPOSITION}
\newtheorem{lema}[teo] {LEMMA}

\newtheorem{obser}[teo]{REMARK}

\newenvironment{demo}{\noindent {\it\bf Proof.} \rm}

\def\Hom{\mathop{\rm Hom}\nolimits}
\def\qed{\hfill \mbox{$\square$}}

\def\op{\mathop{\rm op}\nolimits}

\def\Maps{\mathop{\rm Maps}\nolimits}
\def\Ob{\mathop{\rm Ob}\nolimits}
\def\Mor{\mathop{\rm Mor}\nolimits}
\def\Tot{\mathop{\rm Tot}\nolimits}
\def\L{ {\rm L} }
\def\J{ {\rm J} }
\def\M{ {\rm M} }
\def\tL{ {\tilde {\rm  L}} }
\def\tE{ {\tilde {\rm  E}} }
\def\oE{ {\overline {\rm  E}} }
\def\H{ {\rm H} }
\def\sH{ {\mathbb H} }
\def\E{ {\rm E} }
\def\B{ {\rm B} }
\def\D{ {\rm D} }
\def\F{ {\rm F} }
\def\f{ {\mathcal F}}
\def\C{{\mathbf C}}
\def\Ab{{\mathbf {Ab}}}
\def \Z{{\mathbb Z}}
\def \Q{{\mathbb Q}}
\def\K{{\mathbf K}}
\def\CAT{{\mathbf {CAT}}}
\def\Sets{{\mathbf {Sets}}}

\begin{document}

\title{Cohomology of the Grothendieck construction}
\author{Teimuraz Pirashvili\footnote{T. Pirashvili;  Razmadze Mathematical
Institute of Georgian Academy of Sciences, Alexidze str. 1
Tbilisi, 0193, Republic of Georgia.  E-mail address: pira@rmi.acnet.ge} \ and Mar\'\i a Julia Redondo \footnote{ M. J. Redondo; Instituto de Matem\'atica, Universidad Nacional del Sur, Av. Alem 1253, (8000) Bah\'\i a Blanca, Argentina. E-mail address:  mredondo@criba.edu.ar} \thanks{The first author was supported by DFG at University of Bielefeld; the second author is a researcher from CONICET, Argentina.}}
\date{}
\maketitle

\begin{abstract}
We consider cohomology of small categories with
coefficients in a natural system in the sense of Baues and
Wirsching.  For any functor $\L:\K \to \CAT$, we construct a spectral sequence abutting to the
cohomology of the Grothendieck construction of $\L$ in terms of the cohomology of $\K$ and of $\L(k)$, for  $k \in \Ob \K$.
\end{abstract}


\section{Introduction}
Given a functor $\L:\K \to \CAT$, there is a well defined category
$\int_{\K}\L$ called the Grothendieck construction of $\L$. Our
aim is to construct a spectral sequence abutting to the cohomology of
$\int_{\K} \L$ whose $\E_2$ term is given in terms of the
cohomology of $\K$ and of $\L(k)$, for $k \in \Ob \K$.

If one takes the cohomology of small categories with constant
coefficients then the existence of such spectral sequence is a
consequence of Thomason's work \cite{T}.

In this paper we consider cohomology of small categories with
coefficients in a natural system in the sense of Baues and
Wirsching \cite{BW}.  For a natural system $\D$ on $\int_{\K} \L$
and an object $k$ of $\K$, there is an induced natural system
$\D_k$ on $\L(k)$.  We show that, under some condition on $\D$
(see Theorem \ref{2}), there is a convergent spectral sequence
$$\E_2^{pq} = \H^p(\K, \sH^q(\L(.),\D_{(.)})) \Rightarrow
\H^{p+q}(\mbox{$\int_{\K} \L$}, \D),$$ where
$\sH^q(\L(.),\D_{(.)})$ is a contravariant functor on $\K$ which
assigns to an object $k$ the group $\H^q(\L(k),\D_k)$.  We deduce
this result from another spectral sequence (see Theorem \ref{1})
which holds for any natural system $\D$ but involves cohomology of
some related categories $\tL(k)$ instead of $\L(k)$.

Some particular cases have been already considered in the
literature.  For instance, if $\K$ is the one object category
associated to a group $G$, and the natural system $\D$ comes from
a bifunctor, the mentioned spectral sequence has been constructed
in \cite{CR}, see Remark \ref{cr}.  Another particular case
appears in \cite{P}, where $\K$ is arbitrary but $\L$ has values
in sets, that is, in discrete categories.

\section{Baues-Wirsching cohomology}
Here we fix notation and we prove some elementary facts which will
be needed later. For a small category $\C$ we let $\f\C$ denote
the {\it category of factorizations} in $\C$ \cite{BW}. Objects of
$\f\C$ are the morphisms in $\C$, and a morphism from $\alpha: x
\to y$ to $\beta: u \to v$ is a pair $(\nu: u \to x , \psi: y \to
v)$ of morphisms in $\C$ such that $\beta = \psi \alpha \nu$, that
is, one has a commutative diagram
\[ \xymatrix {   x  \ar[rr]^{\alpha}  &  & y  \ar[d]^{\psi}      \\
      u    \ar[u]^{\nu}  \ar[rr]^{\beta}  &       &    v      } \]
The composition in $\f\C$ is defined by $(\nu, \psi)(\nu',\psi')=
(\nu' \nu, \psi \psi')$. A {\it natural system} on $\C$ is just a
covariant functor $\D: \f\C \to \Ab$.  Now, following \cite{BW},
one defines the cohomology $\H^*(\C,\D)$ as the cohomology of the
cochain complex $\F^*(\C,\D)$ given by
$$\F^n(\C,\D)= \prod_{ c_0  \stackrel{\alpha_1}{\leftarrow}  \cdots
 \stackrel{\alpha_n}{\leftarrow}  c_n} \D( \alpha_1 \cdots \alpha_n) $$
and the coboundary map
$$ d : \F^n(\C,\D) \to \F^{n+1}(\C,\D)$$
is given by
\[
\begin{array}{ll}
(df)(\alpha_1, \cdots, \alpha_{n+1}) &= (\alpha_1)_* f(\alpha_2,
\cdots, \alpha_{n+1})  \\ \\
& + \sum_{i=1}^n (-1)^i f( \alpha_1, \cdots, \alpha_i
\alpha_{i+1}, \cdots, \alpha_{n+1}) \\ \\
& + (-1)^{n+1} (\alpha_{n+1})^* f (\alpha_1, \cdots, \alpha_{n}).
\end{array} \]
Here, and in the rest of the paper, we use the following convention.  For $\alpha
: x \to y$ and elements $a \in \D(u,x), b \in \D(y,v)$, we write
$\alpha_* a$ and $\alpha^*b$ for the image of the elements $a$ and
$b$ by the homomorphisms $\D(id_u, \alpha): \D(u,x) \to \D(u,y)$
and $\D(\alpha, id_v): \D(y,v) \to \D(x,v)$ respectively.

Consider the natural functors
\[ \xymatrix {   \f\C \ar[r] &  \C^{\op} \times \C \ar[r]^{\ \ p_2}  \ar[d]_{p_1}  & \C    \ar[d]   \\
       &     \C^{\op} \ar[r] &   {\bf 1 }    } \]
where the first top functor sends an arrow $\alpha:c \to d$ to the
pair $(c,d)$, $p_1$ and $p_2$ are projections and, finally, $\bf
1$ is the category with one object and one arrow.  Then one gets
that any bifunctor, covariant, contravariant functor or any
abelian group gives rise to a natural system, and then
Baues-Wirsching cohomology generalizes known theories.

Natural systems coming from the functor $\bf 1 \to \Ab$ are called
constant.  In this case, the cohomology of $\C$ is the same as the
cohomology of the corresponding classifying space.

We will need the following easy and well-known result.

\begin{lema}\label{trivial}
Assume $\C$ has an initial object and $\D$ is constant.  Then
\[ \H^n(\C, \D) = \left \{
\begin{array}{ll}
0,  \ & \mbox{if $n>0$},  \\
\D, & \mbox{if $n=0$}.
\end{array} \right. \]
\end{lema}

\begin{demo}
It follows from the fact that the classifying space $\B\C$ is
contractible, see \cite[\S 1, Corollary 2]{Q}. \qed
\end{demo}

\bigskip

The following result is a specialization of \cite[Lemma 1.5, page
10]{FFPS}, setting $\C = \f {\cal A}$, ${\mathbf D}=\f {\cal B}$
and $\F$ the constant natural system on $\C$ with value $\mathbb
Z$.

\begin{lema}\label{adjuntos} Let $({\it l}, {\it r})$ be an
adjoint pair from ${\cal A}$ to ${\cal B}$.  Then for any natural
system $G$ on ${\cal B}$ one has the natural isomorphism
$$\H^*({\cal A}, {\it l}^* (G)) \simeq \H^*({\cal B}, G).$$
\end{lema}

\begin{demo}
Observe first that the induced pair $(\f(\it l), \f(\it r))$ from
$\f{\cal A}$ to $\f {\cal B}$ is also an adjoint pair.  Then we
can use formal adjunction.  \qed
\end{demo}

\bigskip

However we need more sophisticated relationship between
cohomologies of categories involved in adjoint situation.  The
following result is a bit more general than one  proved recently
by Muro (see  \cite[Theorem 6.5]{M}); his proof uses explicit
homotopies,  which he constructs using $2$-categories. Here we use
completely different ideas.

\begin{prop} \label{muro} Let $(\it l, \it r)$  be an adjoint pair from
$\C$ to $\C'$. Assume $\E$ is a natural system on $\C$, and let
$\tE$ and $\oE$ be the natural systems, on $\C$ and $\C'$
respectively, defined by
\[
\begin{array}{lll}
\tE ( c \stackrel{\alpha}{\rightarrow} d) & : = & \E (\epsilon_d \circ \alpha), \\
\oE ( u \stackrel{\beta}{\rightarrow} v) & := & \E (\it r(u)
\stackrel{\it r(\beta)}{\longrightarrow} \it r(v)),
\end{array}
\]
where $\epsilon: id_{\C} \to \it r \it l$ is the unit of the
adjunction.  Then
$$\H^*(\C', \oE) \simeq \H^*(\C, \tE).$$
\end{prop}

\begin{demo}
Varying $\E$ in the category of all natural systems defined on
$\C$, one obtains two sequences of functors
$$(\E \mapsto \H^n(\C, \tE))_{n \geq 0}$$
and
$$(\E \mapsto \H^n(\C, \oE))_{n \geq 0}$$
such that for any short exact sequence of natural systems on $\C$
$$0 \to \E_1 \to \E \to \E_2 \to 0$$
there are corresponding long exact cohomological sequences in both
theories.  Moreover cohomology respects products, therefore it
suffices to show that both theories vanish on injective
cogenerators in positive dimensions and have same values in
dimension zero.

\noindent Let us recall that in a category of functors, the
injective cogenerators are "duals" of representable functors.
Therefore in the category of natural systems on $\C$, the
injective cogenerators are natural systems $\J_{\mu}$, where $\mu:
a \to b$ runs over all morphisms in $\C$ and, for any morphism
$\alpha: c \to d$, one has
$$\J_{\mu}(\alpha)= \Maps (S_{\mu}(\alpha), \Q/\Z)$$
where $S_{\mu}: \f \C ^{\op} \to \Sets$ is given by
$$S_{\mu}(\alpha)= \Hom_{\f\C}(\alpha, \mu)= \left \{ a
\stackrel{\nu}{\longrightarrow} c,  d
\stackrel{\psi}{\longrightarrow} b \mid \psi\alpha\nu=\mu \right
\}.$$ One has an obvious isomorphism of cochain complexes
$$\F^*(\C, \tilde \J_{\mu}) \simeq \F^*(\C_{\mu}, \Q/\Z)$$
where on the right hand side the cohomology is taken with constant
coefficients, while  $\C_{\mu}$ is the category whose objects are
pairs
$$(a \stackrel{\nu}{\longrightarrow} c, \it r \it l ( c)
\stackrel{\psi}{\longrightarrow} b)$$ such that $\psi \circ
\epsilon_c \circ \nu=\mu$, while a morphism from $(\nu, \psi)$ to
$(\nu', \psi')$ is a morphism $\beta: c \to c'$ such that the
following diagrams
\[ \xymatrix
{ & c \ar[dd]^{\beta} & \qquad \qquad & \it r \it l (c)
\ar[dr]^{\psi} \ar[dd]_{\it r \it l (\beta)}
 \\
a \ar[ur]^{\nu} \ar[dr]_{\nu'}  & & \qquad \qquad & & b      \\
& c' & \qquad \qquad &  \it r \it l (c')    \ar[ur]_{\psi'} }
\] commute.  Let us observe also that $\it r \it l (\beta) \circ
\epsilon_c = \epsilon_{c'} \circ \beta$, since $\epsilon $ is a
natural transformation.  Now let $\Delta_{\mu}$ be the set of
morphisms $\gamma: \it r \it l ( a) \to b$ such that $\gamma \circ
\epsilon_a = \mu$.  For any $\gamma \in \Delta_{\mu}$ we let
$\C_{\mu,\gamma}$ be the full subcategory of $\C_{\mu}$ consisting
of objects $(\nu, \psi)$ such that $\psi \circ \it r \it l (\nu) =
\gamma$. It is clear that $\C_{\mu, \gamma} \bigcap \C_{\mu,
\gamma'} = \emptyset$ provided $\gamma \not = \gamma'$.  It is
also clear that for any object $(\nu, \psi) \in \C_{\mu}$ one has
that $(\nu, \psi) \in \C_{\mu,\gamma}$ for $\gamma= \psi \circ \it
r \it l (\nu)$.  Thus $\C_{\mu}$ is the disjoint union
$$\C_{\mu}=\bigsqcup_{\gamma \in \Delta_{\mu}} \C_{\mu, \gamma}.$$
Since for each $\gamma \in \Delta_{\mu}$, $(id_a, \gamma)$ is an
initial object in $\C_{\mu, \gamma}$, we obtain
\[ \H^n(\C, \tilde \J_{\mu}) = \left\{%
\begin{array}{ll}
0, & \hbox{$n>0$;} \\
\Z [\Delta_{\mu}], & \hbox{$n=0$}.
\end{array}%
\right. \] Next we consider the cohomology $\H^*(\C', \overline
\J_{\mu})$. Since $$\F^p(\C', \overline \J_{\mu})= \prod_{u_0
\stackrel{\omega_1}{\leftarrow} \cdots
\stackrel{\omega_p}{\leftarrow} u_p} \overline \J_{\mu}(\omega_1
\circ \cdots \circ \omega_p)$$ and $\overline \J_{\mu}(\omega_1 \circ \cdots
\circ \omega_p)= \prod \Q/\Z$, where the product is taken over all
morphisms $\it r (\omega_1 \circ \cdots \circ
\omega_p) \to \mu$ in $\f \C$, it follows that one has an isomorphism of
complexes
$$\F^*(\C', \overline \J_{\mu}) \simeq \F^*(\C_{\mu}', \Q/\Z)$$
where on the right hand side the cohomology is taken with constant
coefficients and $\C_{\mu}'$ is the category whose objects are
triples
$$(u, a \stackrel{\nu}{\longrightarrow} \it r(u), \it r (u)
\stackrel{\psi}{\longrightarrow} b)$$ such that $\psi \circ
 \nu=\mu$, while a morphism from $(u, \nu, \psi)$ to
$(u',\nu', \psi')$ is a morphism $\tau: u \to u'$ such that the
following diagram
\[ \xymatrix
{ & \it r (u)  \ar[dd]^{\it r(\tau)} \ar[dr]^{\psi}  \\
a \ar[ur]^{\nu} \ar[dr]_{\nu'}  & &   b      \\
& \it r(u') \ar[ur]_{\psi'} }
\] commutes. By adjointness, for any $\nu: a \to \it r (u)$ there exists a
unique $\hat \nu: \it l (a) \to u$ such that the diagram
\[ \xymatrix
{ a    \ar[r]^{\epsilon_a} \ar[dr]_{\nu}  & \it r \it l (a) \ar[d]^{\it r(\hat \nu)}\\
& \it r(u)}
\] commutes. Then it follows that $\C_{\mu}'$ has the following equivalent
description: objects are triples
$$(u, \it l(a) \stackrel{\hat \nu}{\longrightarrow} u, \it r (u)
\stackrel{\psi}{\longrightarrow} b)$$ such that $\mu=\psi \circ
{\it r}(\hat \nu) \circ \epsilon_a$. Recall that $\Delta_{\mu}$ is
the set of morphisms $\gamma: \it r \it l ( a) \to b$ such that
$\gamma \circ \epsilon_a = \mu$.  So, for any $\gamma \in
\Delta_{\mu}$ we let $\C_{\mu,\gamma}'$ be the full subcategory of
$\C_{\mu}'$ consisting of objects $(u, \hat \nu, \psi)$ such that
$\psi \circ \it r (\hat \nu) = \gamma$. It is clear that $\C_{\mu,
\gamma}' \bigcap \C_{\mu, \gamma'}' = \emptyset$ provided $\gamma
\not = \gamma'$.  It is also clear that for any object $(u,\hat
\nu, \psi) \in \C_{\mu}'$ one has that $(u, \hat \nu, \psi) \in
\C_{\mu,\gamma}'$ for $\gamma= \psi \circ \it r (\hat \nu)$. Thus
$\C_{\mu}'$ is the disjoint union
$$\C_{\mu}'=\bigsqcup_{\gamma \in \Delta_{\mu}} \C_{\mu, \gamma}'.$$
Since for each $\gamma \in \Delta_{\mu}$, $(\it l(a), id_{\it
l(a)}, \gamma)= (\it l(a), \hat \epsilon_a, \gamma)$ is an initial
object in $\C_{\mu, \gamma}'$, we obtain
\[ \H^n(\C', \overline \J_{\mu}) = \left\{%
\begin{array}{ll}
0, & \hbox{$n>0$;} \\
\Z [\Delta_{\mu}], & \hbox{$n=0$}.
\end{array}%
\right. \] \qed
\end{demo}

\bigskip

We encountered the following lemma, which we think is of independent interest, when looking for the proof of the previous proposition.

Let $a$ be an object of  a small category $\C$ and let ${\rm
T}:\C^{op}\to \Sets$ be a functor. For any $m\in {\rm T}(a)$ we
define a functor
$$S_{a,{\rm T},m}:\f\C^{op}\to \Sets$$
using the following rule:
$$S_{a,{\rm T},m}(\alpha:c\to d)=\{( a
\stackrel{\eta}{\longrightarrow} c;r\in {\rm T}(d))\mid
\eta^*\alpha^*(r)=m\}.$$ For any abelian group $A$ we put
$$\D_{a,{\rm T},m,A} (\alpha:c\to d)=\Maps(S_{a,{\rm T},m}(\alpha),A).$$
Then $\D_{a,{\rm T},m,A}:\f\C\to \Ab$ is a natural system.

\begin{lema}\label{4vanish} For any small category $\C$, any functor
${\rm T}:\C^{op}\to \Sets$, any element $m\in {\rm T}(a)$ and any
abelian group $A$, one has
\[ \H^n(\C,\D_{a,{\rm T},m,A}) = \left\{%
\begin{array}{ll}
0, & \hbox{$n>0$;} \\
A, & \hbox{$n=0$}.
\end{array}%
\right. \]
\end{lema}

\begin{demo} Consider the category $\C_{a,{\rm T},m}$ whose objects are pairs
$(\eta:a\to d, r\in {\rm T}(d))$, with $\eta^*(r)=m$. A morphism
from $(\eta,r)$ to $(\eta',r')$ in  $\C_{a,{\rm T},m}$ is a
morphism $\beta:d\to d'$ such that $\beta\eta=\eta'$ and
$\beta^*(r')=r$. Then one has an isomorphism of complexes
$$\F^*(\C, \D_{a,{\rm T},m,A})\cong \F^*(\C_{a,{\rm T},m},A).$$
Since $(id_a,m)$ is an initial object of $\C_{a,{\rm T},m}$ the
result  follows. \qed
\end{demo}

\section{Grothendieck construction and Thomason's functor}
Now we recall the classical construction due to Grothendieck
\cite[Expos\'e VI.8]{SGA1}.  Let $\CAT$ be the category of all
small categories.  For any small category $\K$, and any functor
$$\L: \K \to \CAT,$$
the category $\int_{\K}\L$ is defined as follows: objects are
pairs $(k,x)$ with $k \in \Ob \K$ and $x \in \Ob \L(k)$.  A
morphism $(k_0, x_0) \to (k_1, x_1)$ in $\int_{\K}\L$ is a pair
$(\alpha, \xi)$ where $\alpha: k_0 \to k_1$ is a morphism in $\K$
and $\xi : \L(\alpha)(x_0) \to x_1$ is a morphism in $\L(k_1)$.
Composition law is given by
$$(\alpha, \xi)(\alpha', \xi')= (\alpha \alpha', \xi \circ
\L(\alpha)(\xi')).$$ Let $\L: \K \to \CAT$ be a functor.  In
\cite[Definition 1.2.2]{T}, Thomason defines a new functor
$$\tL : \K \to \CAT$$
as follows.  For any object $k$ in $\K$, the category $\tL(k)$ has
as objects pairs $(\alpha, x)$ where $\alpha:l \to k$ is a
morphism in $\K$ and $x$ is an object in $\L(l)$.  A morphism
$(\alpha, x) \to (\alpha', x')$ in $\tL(k)$ is a pair $(\beta,
\xi)$ with $\beta: l \to l'$ and $\xi : \L(\beta)(x) \to x'$, such
that the diagram
\[ \xymatrix {  l  \ar[r]^{\alpha} \ar[d]_{\beta}  &   k       \\
     l' \ar[ur]_{\alpha'} &       } \]
commutes.  Composition is given by
$$(\beta, \xi)(\beta', \xi') = (\beta \beta',
\xi \circ \L(\beta)(\xi')).$$ It is clear that $\tL(k)$ is a well
defined category. Moreover, if $\gamma: k \to k'$ is a morphism in
$\K$, then we have a functor
$$\tL (\gamma) : \tL(k) \to \tL(k')$$
given by
\[
\begin{array}{lll}
\tL (\gamma)(\alpha, x )
& = & (\gamma \alpha, x), \\
\tL(\gamma)(\beta, \xi) & = & (\beta, \xi)
\end{array} \]
for any $(\alpha, x) \in \Ob \tL(k), (\beta, \xi) \in \Mor
\tL(k)$. In this way we get a functor $\tL : \K \to \CAT$.

\section{The first spectral sequence}
Keeping the notation of the preceding section, consider for each
object $k$ in $\K$, the forgetful functor
$$i_k : \tL (k) \to \mbox{$\int_{\K}\L$}$$
given on objects by
$$i_k(\alpha: l \to k, x) = (l,x)$$
and on morphisms by
$$ i_k(\beta, \xi) = (\beta, \xi).$$
If $\D$ is a natural system on $\int_{\K}\L$, the functor $i_k$
induces a natural system $i_k^*\D$ on $\tL (k)$, for any $k \in
\Ob \K$.  Thus we can consider the cohomology
$\H^*(\tL(k),i_k^*\D).$ For any morphism $\gamma: k \to k'$ in
$\K$ one has a commutative diagram of functors
\[ \xymatrix { \tL(k)   \ar[r]^{i_k} \ar[d]_{\tL(\gamma)}  &  \int_{\K}\L      \\
     \tL(k') \ar[ur]_{i_{k'}} &       } \]
In particular, $i_k^* \D = \tL(\gamma)^* i_{k'}^* \D$, and
therefore we have a well defined homomorphism in cohomology
$$\H^*(\tL(k),i_k^*\D) \to \H^*(\tL(k'),i_{k'}^*\D).$$
In fact, we have a functor,
$$\sH^*(\tL(.),i_{(.)}^*\D) : \K^{\op} \to \Ab,$$
sending $k \in \Ob (\K)$ to $\H^*(\tL(k),i_k^*\D)$.

\begin{teo}\label{1}
Let $\K$ be a small category, $\L : \K \to \CAT$ a functor and
$\D$ a natural system  on $\int_{\K}\L $.  Then there is a first
quadrant spectral sequence
$$ \E_2^{pq} = \H^p(\K, \sH^q(\tL(.), i_{(.)}^* \D) \Rightarrow
\H^{p+q}(\mbox{$\int_{\K}\L$} , \D).$$
\end{teo}

\begin{demo}
The proof is inspired by an argument given by Thomason in
\cite[Lemma 1.2.4]{T}.  Let $C^{**}$ be the bicomplex defined as
follows.  For any $p,q \geq 0$ let $M^{p,q}$ be the set of all
triples
$$( i_0 \stackrel{\alpha_1}{\leftarrow} \cdots
\stackrel{\alpha_p}{\leftarrow}i_p \  ; \  (j_0, x_0)
\stackrel{(\beta_1, \xi_1)}{\longleftarrow} \cdots
\stackrel{(\beta_q, \xi_q)}{\longleftarrow}(j_q,x_q) \ ; \ i_p
\stackrel{\gamma}{\leftarrow} j_0)$$ consisting of $p$ composable
morphisms $\alpha_1, \cdots, \alpha_p$ in $\K$, $q$ composable
morphisms $(\beta_1, \xi_1), \cdots , (\beta_q, \xi_q)$ in
$\int_{\K}\L$, and a morphism in $\K$ connecting the first
component of the target of $(\beta_1, \xi_1)$ with the source of
$\alpha_p$.  Now we set
$$ C^{p,q} = \prod \D ((\beta_1, \xi_1) \circ \cdots \circ (\beta_q, \xi_q))$$
where the product is taken over the set $M^{p,q}$. The
homomorphisms
$$ C^{p,q}  \stackrel{\delta^{p,q}}{\longrightarrow}  C^{p+1,q}  \qquad \mbox{and}
\qquad C^{p,q}  \stackrel{\partial^{p,q}}{\longrightarrow}  C^{p,
q+1}$$ are given by
{\small
\[
\begin{array}{ll}
& \delta^{p,q}(f)  (\alpha_1, \cdots, \alpha_{p+1} \  ; \
(\beta_1,
\xi_1), \cdots , (\beta_q, \xi_q) \ ; \ \gamma ) =  \\ \\
& = f(\alpha_2, \cdots, \alpha_{p+1} \ ; \ (\beta_1,
\xi_1), \cdots , (\beta_q, \xi_q) \ ; \ \gamma ) \\ \\
& + \sum_{t=1}^p (-1)^t f( \alpha_1, \cdots, \alpha_t
\alpha_{t+1}, \cdots, \alpha_{p+1} \ ; \ (\beta_1,
\xi_1), \cdots , (\beta_q, \xi_q) \ ; \ \gamma ) \\ \\
& + (-1)^{p+1}  f (\alpha_1, \cdots, \alpha_{p} \ ; \ (\beta_1,
\xi_1), \cdots , (\beta_q, \xi_q) \ ; \ \alpha_{p+1} \gamma )
\end{array} \]}
and {\small \[
\begin{array}{ll}
& \partial^{p,q}(f) (\alpha_1, \cdots, \alpha_{p} \ ; \ (\beta_1,
\xi_1), \cdots , (\beta_{q+1}, \xi_{q+1}) \ ; \ \gamma ) = \\ \\&
= (\beta_1, \xi_1)_* f(\alpha_1, \cdots, \alpha_{p} \ ; \
(\beta_2,
\xi_2), \cdots , (\beta_{q+1}, \xi_{q+1}) \ ; \ \gamma \beta_1) \\ \\
& + \sum_{t=1}^q (-1)^t f( \alpha_1, \cdots, \alpha_{p} \ ; \
(\beta_1,
\xi_1), \cdots ,(\beta_t, \xi_t)(\beta_{t+1}, \xi_{t+1}), \cdots,  (\beta_{q+1}, \xi_{q+1}) \ ; \ \gamma ) \\ \\
& + (-1)^{q+1} (\beta_{q+1}, \xi_{q+1})^* f (\alpha_1, \cdots,
\alpha_{p} \ ; \ (\beta_1, \xi_1), \cdots , (\beta_q, \xi_q) \ ; \
\gamma ).
\end{array} \]}
Let $\Tot(C^{**})$ be the total complex associated to the
bicomplex $C^{**}=(C^{p,q}, \delta^{p,q}, (-1)^p \partial^{p,q})$.
First we define the morphism of complexes
$$\phi^* : \F^*(\mbox{$\int_{\K}\L$}, \D) \to C^{0,*}$$
using the following formula
$$\phi^q  (f)( i_0 \  ; \  (j_0, x_0)
\stackrel{(\beta_1, \xi_1)}{\longleftarrow} \cdots
\stackrel{(\beta_q, \xi_q)}{\longleftarrow}(j_q,x_q) \ ; \ i_0
\stackrel{\gamma}{\leftarrow} j_0) = f( (\beta_1, \xi_1) \circ
\cdots \circ (\beta_q, \xi_q)).$$ Since the composition
$$\F^*(\mbox{$\int_{\K}\L$}, \D)
\stackrel{\phi^*}{\longrightarrow} C^{0,*}
\stackrel{\delta^{0,*}}{\longrightarrow} C^{1,*}$$ vanishes,
$\phi^*$ induces a homomorphism of complexes $\overline {\phi^*}
: \F^*(\mbox{$\int_{\K}\L$}, \D) \to \Tot^*(C^{**})$. We claim
that $\overline {\phi^*}$ is a quasi-isomorphism.  Indeed, it
suffices to show that for each $n \geq 0$ the sequence
$$ 0 \to \F^n(\mbox{$\int_{\K}\L$}, \D) \to C^{0,n}
\stackrel{\delta^{0,n}}{\longrightarrow} C^{1,n}
\stackrel{\delta^{1,n}}{\longrightarrow} C^{2,n}
\stackrel{\delta^{2,n}}{\longrightarrow}  \cdots $$ is exact. Note
that the complex $C^{*,n}$ can be represented as the product over
all $n$-strings
$$ (j_0, x_0)
\stackrel{(\beta_1, \xi_1)}{\longleftarrow} \cdots
\stackrel{(\beta_n, \xi_n)}{\longleftarrow}(j_n,x_n)$$ of the
complexes
$$\F^*(\K \uparrow j_0, D(  (\beta_1, \xi_1) \circ \cdots \circ
(\beta_n, \xi_n))),$$ where $\K \uparrow j_0$ is the category
whose objects are morphisms $j \stackrel{\sigma}{\leftarrow} j_0$
in $\K$ and morphisms $\sigma \to \sigma'$ are morphisms $\tau: j
\to j'$ in $\K$ such that the following diagram commutes
\[ \xymatrix {  j \ar[d]_{\tau}  &  j_0  \ar[l]_{\sigma} \ar[dl]^{\sigma'}   \\
        j' &        } \]
Let us observe that the coefficient $D(  (\beta_1, \xi_1) \cdots
(\beta_n, \xi_n))$ is constant and $\K \uparrow j_0$ has an
initial object $id_{j_0}$, so the claim follows from Lemma
\ref{trivial}. Now, it follows from the spectral sequence of a
bicomplex that to finish the proof we only need to prove that
$$\H^*(C^{p,*})= \F^p(\K, \sH^*(\tL(.), i_{(.)}^* \D ).$$
But $C^{p,*}$ is the product over all $p$-strings $i_0
\stackrel{\alpha_1}{\leftarrow} \cdots
\stackrel{\alpha_p}{\leftarrow}i_p$ of the complexes
$\F^*(\tL(i_p), i^*_{i_p}\D)$, and the result follows. \qed
\end{demo}

\section{The second spectral sequence}
Let $\K$ be a small category and let $\L: \K \to \CAT$ be a
functor as above. For each object $k \in \Ob \K$, one has an
adjoint pair $({\it l}_k, {\it r}_k)$ from $\tL(k)$ to $\L(k)$
given by
$$ {\it l}_k (\alpha, x) =  \L(\alpha)(x) \qquad \mbox{and} \qquad
{\it r}_k(y)  =  (id_k, y)$$ for any $(l
\stackrel{\alpha}{\rightarrow} k, x \in \Ob \L(l)) \in \Ob \tL(k)$
and any $y \in \Ob \L(k)$. For a natural system $\D$ on
$\int_{\K}\L$ and an object $k$ of $\K$, we have the induced
natural system $i^*_k\D$ on  $\tL(k)$. Since $({\it l}_k, {\it
r}_k)$ is an adjoint pair from $\tL(k)$ to $\L(k)$, we can apply
Proposition \ref{muro}, setting $\E_k= i_k^*\D$. One has $\oE_k =
\D_k$, where $\D_k$ is a natural system on $\L(k)$ given by
$$\D_k( x_0 \stackrel{\xi}{\rightarrow} x_1) =
\D((k,x_0) \stackrel{(id_k, \xi)}{\longrightarrow}(k,x_1)).$$ On
the other hand
$$\tE_k(\beta, \psi)= \D(\alpha, \L(\alpha')(\psi))$$
for any morphism $(\beta, \psi): (\alpha,x) \to (\alpha',x')$ in
$\tL(k)$.  Hence Proposition \ref{muro} implies
\begin{equation}\label{kadjiso}
\H^*(\L(k), \D_k) \simeq \H^*(\tL(k), \tE_k ).
\end{equation}
On the other hand there is a natural transformation $i^*_k \D \to
\tE_k $ of natural systems on $\tL(k)$. We will say that the
natural system $\D$ is \emph{h-local} provided, for any object $k$
in $\K$, the induced homomorphism in cohomology $\H^*(\tL(k),i^*_k
\D)\to  \H^*(\tL(k),\tE_k)$ is an isomorphism. If this holds then
$\H^*(\tL(k),i^*_k \D)\cong \H^*(\L(k), \D_k)$ thanks to
isomorphism (\ref{kadjiso}). Thus in this case the functors
$\sH^*(\L(.), \D_{( . )})$ and $\sH^*(\tL(.), i^*_{( . )} \D)$ are
isomorphic. By  Theorem \ref{1} we obtain the following result.

\begin{teo}\label{2}
Let $\K$ be a small category and let $\L: \K \to \CAT$ be a
functor.  Assume $\D$ is a natural system on $\int_{\K}\L$ which
is h-local. Then there is a first quadrant spectral sequence
$$\E_2^{pq} = \H^p(\K, \sH^q(\L(.),\D_{(.)})) \Rightarrow
\H^{p+q}(\mbox{$\int_{\K} \L$}, \D),$$ where
$\sH^q(\L(.),\D_{(.)}): \K^{\op} \to \Ab$ sends $k \in \Ob
\K$ to $\H^q(\L(k),\D_{k})$ and $\D_k$ is a natural system on
$\L(k)$ given by $$\D_k( x_0 \stackrel{\xi}{\rightarrow} x_1) =
\D((k,x_0) \stackrel{(id_k, \xi)}{\longrightarrow}(k,x_1)).$$
\end{teo} \qed

\begin{obser} \rm Assume $\D$ is a natural system such that the morphism
$i^*_k \D \to \tE_k $ of natural systems is an isomorphism. Then
obviously $\D$ is h-local. In this case we will say that $\D$ is
\emph{local}.  Then $\D$ is local if and only if for any morphism
$(k_0, x_0)\stackrel{(\alpha, \xi)}{\longrightarrow} (k_1, x_1)$
in $\int_{\K}\L$ and any arrow $\beta: k_0 \to k_1$, the induced
homomorphism
$$(\beta, id_{\L(\beta)(k_1)})_* : \D(\alpha, \xi) \to \D((\beta,
id_{\L(\beta)(k_1)})_*(\alpha, \xi))$$ is an isomorphism. It is
now clear that this condition automatically holds, for instance,
when $\K$ is a groupoid, or $\K$ is arbitrary and $\D$ is of the
form
$$\D((k_0,x_0) \stackrel{(\alpha,
\xi)}{\longrightarrow}(k_1,x_1))= \M(k_0,x_0)$$ where $\M:
(\int_{\K} \L)^{\op} \to \Ab$ is a functor.
\end{obser}

\begin{obser}\label{cr} \rm
Let $k$ be a field and let $\cal C$ be a $k$-linear category. If
$G$ is a group acting freely on ${\cal C}$, one can consider the
quotient category ${\cal C} / G$.  Now, for any locally finite
${\cal C} / G$-bimodule  $\cal M$, there is a spectral sequence
$$ \E_2^{p,q} = H^p (G, H^q({\cal C}, \L {\cal M}))
\Rightarrow H^{p+q}({\cal C} / G , {\cal M}),$$ see \cite[Theorem
3.11]{CR}. Moreover, the skew category ${\cal C}[G]$ defined in
\cite{CM} is equivalent to the quotient category ${\cal C}/ G$,
see \cite[Theorem 2.8]{CM}. Assume now that $\K$ is the one object
category associated to the group $G$.  Define $\L$ as the functor
$\K \to \CAT$ that sends the unique object to the category $\cal
C$, and morphisms in $\K$ corresponding to the $G$-action. One easily
sees that for the corresponding $k$-linear span $k[{\cal C}]$ one
has an isomorphism of $k$-linear categories
$$k[{\cal C}][G] \simeq k[ \mbox{$\int_{\K}\L$}].$$
Thus the spectral sequence constructed in \cite{CR} gives the
spectral sequence constructed in \ref{2} for such $\K$, when
coefficients are bifunctors.
\end{obser}

\begin{obser} \rm
Let $\K$ be the one object category associated to a group $G$ and
let $\L : \K \to \CAT$ be a functor, that is, a category $\C$ with
a $G$-action.  For any natural system $\D$ on $\int_{\K}\L$ ,
Theorem \ref{2} says that we have a spectral sequence
$$ \E_2^{p,q} = H^p (G, H^q(\C, \D))
\Rightarrow H^{p+q}(\mbox{ $\int_{\K}\L$},\D).$$ In this
particular case, an alternative proof can be given using a
well-known theorem of Eilenberg-Zilber-Cartier \cite{DP} that says
that the diagonal of a bicosimplicial abelian group is homotopy
equivalent to the total complex of the associated bicomplex.
Consider the bicomplex $C^{**}= C^*(G, F^*(\C, \D))$, where
$C^*(G,M)$ denotes the Eilenberg-Mac Lane cochains of the group
$G$ with coefficients in a $G$-module $M$, and the $G$-action on
the cochain complex $\F^*(\C, \D)$ is given by
$$f^s(\alpha_1, \cdots, \alpha_p)=
(s,id)_*(s^{-1},id)^*f(^{s^{-1}} \alpha_1, \cdots, ^{s^{-1}}
\alpha_p)$$ for any $s \in G$. An explicit isomorphism from the
diagonal of $C^{**}$ to $\F^*(\int_{\K}\L, \D)$ is given as
follows. For any $g_1, \cdots, g_p \in G$, $\alpha_1, \cdots,
\alpha_p \in \C$, we set
$$\phi^p(f)((g_1, \alpha_1) \circ \dots \circ
(g_p, \alpha_p)) =$$ $$=(g_1 \cdots g_p, id)^* f (g_1, \dots, g_p,
\alpha_1,  ^{g_1} \alpha_2, \cdots, ^{g_1 \cdots g_{p-1}}
\alpha_p).$$
\end{obser}

\bigskip

\noindent {\bf ACKNOWLEDGEMENTS.} This paper was written while the
authors were visiting the Universit\'e de Montpellier II.  They
want to express their gratitude to Claude Cibils and Daniel Guin
for their hospitality.

\end{document}